\documentclass[12pt]{article}
\usepackage[numbers]{natbib}
\usepackage{graphicx}
\usepackage{color, multirow,amsfonts,amsmath}
\usepackage{url}
\usepackage{rotating}
\title{Mixed Poisson process with Max-U-Exp mixing variable - Working version
\bigskip
\\Pavlina K. Jordanova
\\{\small{\it Faculty of Mathematics and Informatics, Konstantin Preslavsky University of Shumen, \\115 "Universitetska" str., 9712 Shumen, Bulgaria. \\Corresponding author:  pavlina\_kj@abv.bg.}}
\\Evelina Veleva
\\{\small{\it Department of Applied mathematics and Statistics, "Angel Kanchev" University of Ruse, Bulgaria.}}
}

\begin{document}
\date{}

\maketitle

\begin{abstract}
This work defines and investigates the properties of the Max-U-Exp distribution. The method of moments is applied in order to estimate its parameters. Then, by using the previous general theory about Mixed Poisson processes, developed by Grandel (1997), and Karlis and Xekalaki (2005), and analogously to Jordanova et al. (2023), and Jordanova and Stehlik (2017) we define and investigate the properties of the new random vectors and random variables, which are related with this particular case of a Mixed Poisson process.  Exp-Max-U-Exp distribution is defined and thoroughly investigated. It arises in a natural way as a distribution of the inter-arrival times in the Mixed Poisson process with Max-U-Exp mixing variable. The distribution of the moments of arrival of different events is called Erlang-Max-U-Exp and is intorduced via its probability density function. Investigation of its properties follows. Finally, the corresponding Mixed Poisson process with Max-U-Exp mixing variable is defined. Its finite dimensional and conditional distributions are found and their numerical characteristics are determined.
\end{abstract}

\section{INTRODUCTION}

The total set of probability distributions and random processes is uncountable, therefore, when introduce and investigate them it is desirable to show the connections between them. One way to make this, is to start with some random process, and to obtain all random distributions of the stochastic elements which describe it. In 1997 Grandel \cite{GrandelMixed} summarised and developed the general theory of Mixed Poisson processes and present some of their potential applications. Later on, in 2005, Karlis and Xekalaki\cite{KarlisandXekalaki} make a very good review of the investigations of many particular cases of such processes and obtain some of their new properties and multivariate versions. Analogously, in 2017 Jordanova and Stehlik \cite{LMJ} study the case when the mixing variable is Pareto distributed and define the distributions which describe, the univariate and multivariate distributions processes related with this case, the distribution of the inter-arrival times, the moment of the $n$-the event and so forth. In 2023 Jordanova et al. \cite{SMladen} consider the very useful and general case when the mixing variable is Stacy distributed. Here we define a new Max-U-Exp distribution and investigate its properties. The method of moments is applied in order to estimate its parameters. Then, by using the previous general theory about Mixed Poisson processes we define and investigate the properties of the new random vectors and random variables, which are related to this particular case of a Mixed Poisson process.  Exp-Max-U-Exp distribution is defined and thoroughly investigated. It arises in a natural way as a distribution of the inter-arrival times in the Mixed Poisson process with Max-U-Exp mixing variable. The distribution of the renewal moments is called Erlang-Max-U-Exp and is defined via its probability density function. Investigation of its properties follows. Finally, the corresponding Mixed Poisson process with Max-U-Exp mixing variable is defined. Its finite dimensional and conditional distributions are found and their numerical characteristics are determined.

Along this work we denote by $\xi \in Bi(n, p)$ the fact that a random variable (r.v.) $\xi$ belongs to the set of Binomial distributions with parameters $n \in \mathbb{N}$, and $p \in (0, 1)$. As usually $\Gamma(\alpha) = \int_0^\infty x^{\alpha-1}e^{-x}dx$ is the notation for the Euler's Gamma function.  $\Gamma(\alpha, x) = \int_x^\infty y^{\alpha-1}e^{-y}dy$, and $\gamma(\alpha,x) = \Gamma(\alpha) - \Gamma(\alpha,x)$ are correspondingly the upper, and the loweer incomplete Gamma functions. $F_\xi(x)$ is the cumulative distribution function (c.d.f.) of the r.v. $\xi$ and  $P_\xi(x)$ is for its probability density function (p.d.f.).

\section{MAX-U-EXP DISTRIBUTION}

{\bf Definition 1.} \label{Def:1} We say that the r.v. $\xi$ is {\bf{Max-U-Exp distributed}} with parameters $a > 0$ and $\lambda > 0$, if it has a cumulative distribution function (c.d.f.)
\begin{equation}\label{MaxUExpCDF}
F_\xi(x) = \left\{ \begin{array}{ccc}
                                   0 & , & x \leq 0\\
                                   \frac{x}{a}(1-e^{-\lambda x})& , & x \in (0, a]\\
                                   1-e^{-\lambda x}& , & x > a\\
                      \end{array}\right..
\end{equation}

Briefly we will denote this in this way  $\xi \in Max-U-Exp(a; \lambda)$.

{\bf Proposition 1.}
\begin{description}
  \item[a)] $\xi \in Max-U-Exp(a; \lambda)$ if and only if the probability density function of $\xi$ is
\begin{equation}\label{MaxUExpDensity}
P_\xi(x) = \left\{ \begin{array}{ccc}
                                   0 & , & x \leq 0\\
                                   \frac{1}{a}(1-e^{-\lambda x}+ x \lambda e^{-\lambda x})& , & x \in (0, a]\\
                                   \lambda e^{-\lambda x}& , & x > a\\
                      \end{array}\right..
                      \end{equation}
  \item[b)] (Scaling property) If $\xi \in Max-U-Exp(a; \lambda)$ and $k > 0$ is a constant, then
  $$k\xi \in Max-U-Exp\left(ka; \frac{\lambda}{k}\right).$$
  \item[c)]  If $\xi \in Max-U-Exp(a; \lambda)$, the hazard rate function of this distribution is
  $$h_{\xi}(x) = \left\{\begin{array}{ccc}
                   0 & , & x \leq 0 \\
                   \frac{1-e^{-\lambda x}+ x \lambda e^{-\lambda x}}{a-x + x e^{-\lambda x}} & , & x \in (0, a]\\
                    \lambda & , & x > a
                  \end{array}
  \right.$$
\end{description}

{Proof:} a) is an immediate corollary of the relation between c.d.f. and probability density function (p.d.f.).

b) For $k > 0$ and $x > 0$, by using a) we obtain
$$ P_{k\xi}(x) =  \frac{1}{k}P_{\xi}\left(\frac{x}{k}\right) = \left\{ \begin{array}{ccc}
                                   0 & , & x \leq 0\\
                                   \frac{1}{ak}(1-e^{-\frac{\lambda}{k} x}+ x \frac{\lambda}{k} e^{-\frac{\lambda}{k} x})& , & x \in (0, ka]\\
                                   \frac{\lambda}{k} e^{-\frac{\lambda}{k} x}& , & x > ka\\
                      \end{array}\right..$$
The rest follows by the uniqueness of the correspondence between p.d.f. and the distribution and formula (\ref{MaxUExpDensity}).

c) follows by the definition for hazard rate function $h_{\xi}(x) = \frac{P_{\xi}(x)}{1-F_{\xi}(x)}$, Definition 1, and Proposition 1, a).
 \hfill $\Box$

In the next theorem and further on, we denote by $U(0, a)$ the Uniform distribution on the interval $(0, a)$, and by $Exp(\lambda)$ the Exponential distribution with mean $\frac{1}{\lambda}$, $\lambda > 0$.

 {\bf Theorem 1.} Let $\theta \in U(0, a)$, $\eta \in Exp(\lambda)$ and $\theta$ and $\eta$ be independent.
 Denote by $\xi := \max(\theta, \eta)$. Then,
\begin{itemize}
\item [a)] $\xi \in Max-U-Exp(a; \lambda)$;
\item [b)] The mean, and the moments of $\xi$ are correspondingly $\mathbb{E}\xi = \frac{a}{2} + \frac{1}{a \lambda^2}(1 - e^{-\lambda a}),$ and
$$\mathbb{E}(\xi^k) = \frac{a^k}{k+1} + \frac{k}{a\lambda^{k+1}}\gamma(k+1, a\lambda) + \frac{k}{\lambda^k}\Gamma(k, \lambda a), \quad k > -1.$$
\item [c)] The variance of $\xi$ is $\mathbb{D}\xi = \frac{a^2}{12} - \frac{1}{\lambda^2}(1 + e^{-\lambda a}) + \frac{4}{a \lambda^3}(1 - e^{-\lambda a})  -\frac{1}{a^2 \lambda^4}(1-e^{-\lambda a})^2 .$
\item [d)] The Laplace-Stieltjes transform of $\xi$ is
  $$\mathbb{E}(e^{-\xi t}) = \frac{1}{at}(1 - e^{-\lambda a}) - \frac{t}{a(\lambda + t)^2}(1 - e^{-(\lambda + t) a}).$$
\end{itemize}

{Proof:} a) Consider $x \in \mathbb{R}$, the definition of $\xi$ and the independence between $\theta$ and $\eta$ entail,
$$F_\xi(x) =  \mathbb{P}(\max(\theta, \eta) \leq x) = \mathbb{P}(\theta \leq x, \eta \leq x)  = \mathbb{P}(\theta \leq x) \mathbb)\mathbb{P}(\eta \leq x)$$
Now by using the definitions of $Exp(\lambda)$ and $U(0, a)$ distributions we obtain the c.d.f. (\ref{MaxUExpCDF}). The rest follows by the uniqueness of the correspondence between the c.d.f. and the probability law of the considered random variable (r.v.).

b) follows by the definition of the mathematical expectation, initial moments,  and (\ref{MaxUExpDensity}).

c) follows by the formula $\mathbb{D}\xi = \mathbb{E}(\xi^2) - (\mathbb{E}\xi)^2$, the definition of the second initial moment of a r.v.,  and (\ref{MaxUExpDensity}).

d) is a corollary of the definition for Laplace-Stieltjes transform of a r.v.,  and (\ref{MaxUExpDensity}).  \hfill $\Box$

\begin{figure}[ht]
  \centering
  \includegraphics[scale=.45]{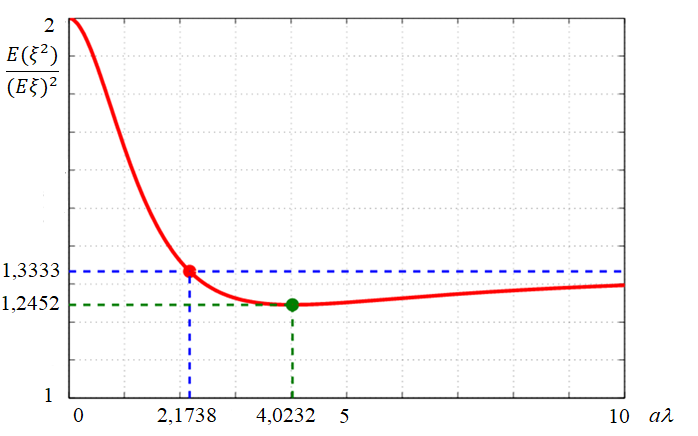}
  \caption{ \label{fig:Fig1Eva}}
\end{figure}

Let us now use the method of moments, and to obtain the algorithm for estimation of the parameters of this distribution. Suppose we have a sample of $n$ independent observations on a r.v. $\xi$.  Let us denote by $m_k$ the $k$-th empirical initial moment of $\xi$ computed by using a sample of $n$ independent observations on a r.v. $\xi \in Max-U-Exp(a; \lambda)$. Then, it is well-known that $m_1$, and $m_2$  are unbiased and consistent estimations correspondingly for $\mathbb{E}\xi$, and $\mathbb{E}(\xi^2)$, while an unbiased estimation for $(\mathbb{E}\xi)^2$ is $\frac{nm_1^2-m_2}{n-1}$. By Theorem 1, we obtain that the first two initial moments are the following functions of the unknown parameters $a$ and $\lambda$:
$$\left|\begin{array}{c}
E\xi = \frac{a}{2} + \frac{1}{a\lambda^{2}}\left(1 - e^{-\lambda a}\right) = \frac{1}{xy}\left(\frac{x^2}{2}+1-e^{-x}\right)\\
E(\xi^2) =  \frac{a^2}{3} + \frac{4}{a\lambda^{3}}(1-e^{-a\lambda}) - \frac{2}{\lambda^2}e^{-\lambda a} = \frac{1}{xy^2}\left[\frac{x^3}{3}+4-{2e}^{-x}(x+2)\right]\\ \end{array}
\right. $$
where we have used the notations $x = a\lambda$, $y = \lambda$. Equivalently,
\begin{equation}\label{3}
\left|\begin{matrix}x\left[\frac{x^3}{3}+4-{2e}^{-x}(x+2)\right]:\left(\frac{x^2}{2}+1-e^{-x}\right)^2 = \frac{E\xi^2}{{(E\xi)}^2}\\y = \left(\frac{x^2}{2}+1-e^{-x}\right):(xE\xi)\\\end{matrix}\right.
\end{equation}
The first equation of system (3) is nonlinear, depending only on the unknown $ x =a \lambda $. Its solution can be found numerically by replacing its right-hand side with an estimate $\hat{r}$ for the ratio $\frac{\mathbb{E}(\xi^2)}{(\mathbb{E}\xi)^2}$ calculated from the sample. Such estimate can be  $\frac{m_2}{m_1^2}$  or $\frac{m_2 (n-1)}{nm_1^2-m_2}$, $\frac{m_2}{m_1^2} > \frac{m_2 (n-1)}{nm_1^2-m_2}$. After finding $x$,  we determine the unknown $y = \lambda$, estimating $\mathbb{E}\xi$ with $m_1$ from the second equation of the system (3). Finally, we find $a = \frac{x}{y}$. The graph of the left-hand side of the first equation of system (3) as a function of $x = a\lambda$ is shown in Figure \ref{fig:Fig1Eva}. We can see that when

the estimator $\hat{r} \in \left[\frac{4}{3}, 2\right]$, the system (3) will have a unique solution. It can easily be checked that as $x = a\lambda$  tends to infinity the ratio $\frac{\mathbb{E}(\xi^2)}{(\mathbb{E}\xi)^2}$ will tend to $\frac{4}{3} = 1.3333$. However, to each $\hat{r}$ in the interval $[1.2452,1.3333]$ we have two possible values of $a\lambda$ on the graph of the left side of the first equation, i.e. two possible solutions. The constant $1.2452$ is the minimum value that the ratio $\frac{\mathbb{E}(\xi^2)}{(\mathbb{E}\xi)^2}$ can take. It is reached when $a\lambda = 4.0232$. In the simulations made, sometimes $\hat{r}$ took values even less than 1.2452. In such a case, the system (3) will not have a solution, or we can assume that $a\lambda = 4.0232$ as the value that minimizes the square of the difference between the left and right sides of the first equation.

When $\hat{r} < 1.3333$ we could conclude that $a\lambda > 2.1738$ (the $x$-coordinate of the corresponding value on the graph) and use another approach to estimate $a\lambda$, for example the least squares method to compare the empirical with the theoretical distribution functions. In detail, this approach is described for instance in \cite{Rameshwar} for parameter estimation in the Generalized exponential distribution. From (1) we have that $\mathbb{P}(\xi > a) = 1 - F_\xi(a) = e^{-\lambda a}$, and for $a\lambda > 2.1738$, $\mathbb{P}(\xi > a)$ will be less than 0.1137, that is, a relatively small percentage of the observations in the sample will be greater than the parameter $a$. With probability close to one no more than 25\% of the observations will be greater than $a$. For example, if $n = 20$ and $a\lambda = 2.2$, $\mathbb{P}(\xi > a) = e^{-2.2} = 0.11$, and the probability that no more than 25\% of the observations are greater than $a$ is
$$P(X \leq 5) = \sum_{i=0}^5 C_{20}^i 0.11^i 0.89^{20-i} = 0.98,$$
where $X \in Bi(20, 0.11)$. For $a\lambda >2.1738$  we can simply remove the largest 25\% of the observations and with the remaining $k$ observations $x_1, x_2, \ldots, x_k$ form and minimize with respect to the parameters $a$ and $\lambda$ the sum
\begin{equation}\label{4}
 \sum_{i=1}^{k}\left[\frac{i}{n+1}-F_\xi(x_i)\right]^2 = \sum_{i=1}^{k}\left[\frac{i}{n+1}-\frac{x_i}{a}\left(1-e^{-\lambda x_i}\right)\right]^2.
\end{equation}

Alternatively, a lower initial estimate $\hat{a}$ for the parameter $a$ can be determined using the histogram of the sample. Since the density function of the distribution will always have a discontinuity at point $a$, from a given location $\hat{a}$ onwards the heights of the bars in the histogram will drop sharply, decreasing exponentially to $0$. For larger values of the product $a\lambda$, the exponential tail will not even be present in the histogram at all, and all observations will be in the interval $(0, \hat{a})$. Then, we determine the number $k$ of elements in the sample that are smaller than $\hat{a}$ and with their help we form and minimize the sum (4) with respect to the parameters $a$ and $\lambda$.
\bigskip

\section{EXP-MAX-U-EXP AND ERLANG-MAX-U-EXP DISTRIBUTIONS}

{\bf Definition 2.} \label{Def:2} We say that the r.v. $\tau$ is {\bf{Exp-Max-U-Exp distributed}} with parameters $a > 0$ and $\lambda > 0$, if it has a p.d.f.
\begin{equation}\label{ExpMaxUExpDensity}
P_\tau(t) = \left\{ \begin{array}{ccc}
                                   0 & , & t \leq 0\\
                                   \frac{1}{at^2}\left(1-e^{-at}-at e^{-at}\right)+\frac{\lambda - t}{a(\lambda + t)^3}(1-e^{-a(\lambda + t)})+\frac{t}{(\lambda+t)^2}e^{-a(\lambda + t)}& , & t > 0
                      \end{array}\right..
\end{equation}

Briefly we will denote this in this way  $\tau \in Exp-Max-U-Exp(a; \lambda)$.

This distribution is proper as far as $\int_0^\infty P_\tau(t) dt = 1$.

The proof of the following result is based on  the correspondence between c.d.f., p.d.f. and the probability distribution.

{\bf Proposition 2.} For $a > 0$ and $\lambda > 0$, $\tau \in Exp-Max-U-Exp(a; \lambda)$ if and only if the c.d.f.
\begin{equation}\label{ExpMaxUExpCDF}
F_\tau(t) = \left\{ \begin{array}{ccc}
                                   0 & , & t \leq 0\\
                                   1-\frac{1-e^{-at}}{at}+\frac{t}{a(\lambda + t)^2}(1-e^{-a(\lambda + t)})& , & t > 0
                      \end{array}\right..
\end{equation}

{\bf Definition 3.} \label{Def:3} We say that the random vector (rv.) $(\tau, \xi)$ has {\bf{bivariate Exp-Max-U-Exp distribution of $I^{-st}$ kind}} with parameters $a > 0$, and $\lambda > 0$, if it has a joint p.d.f.
\begin{equation}\label{BivariateExpMaxUExpDensity}
P_{\tau, \xi}(t,x) = \left\{ \begin{array}{ccc}
                                   0 & , & x \leq 0 \cup t \leq 0\\
                                   \frac{xe^{-tx}}{a}(1-e^{-\lambda x}+ x \lambda e^{-\lambda x})& , & x \in (0, a], t > 0\\
                                   \lambda x e^{-(\lambda + t) x}& , & x > a, t > 0\\
                      \end{array}\right..
\end{equation}
Briefly we will denote this in this way  $(\tau, \xi) \in Exp-Max-U-Exp-I^{-st}(a, \lambda)$.

{\bf Theorem 2.} For $a > 0$ and $\lambda > 0$, if $\xi \in Max-U-Exp(a; \lambda)$ and for $x > 0$, $(\tau|\xi=x) \in Exp(x)$, then:

\begin{description}
  \item[a)] $\tau \in Exp-Max-U-Exp(a, \lambda)$;
  \item[b)] $\tau \stackrel{d}{=} \frac{\eta}{\xi}$, where $\eta \in Exp(1)$, and $\xi$ and $\eta$ are independent.
  \item[c)] For $p \in (0, 1)$,
  $$\mathbb{E}(\tau^p) = \Gamma(p+1)\left(\frac{1}{a^p(1-p)}+\frac{\lambda^{p-1}}{a}\left((p+a\lambda)\Gamma(1-p,\lambda a)+(\lambda a)^{1-p}e^{-\lambda a} -p\Gamma(1-p)\right)\right).$$
  For $p \geq 1$, $\mathbb{E}(\tau^p) = \infty$.
  \item[d)] The joint distribution of $\tau$ and $\xi$ is $(\tau, \xi) \in Exp-Max-U-Exp-I^{-st}(a, \lambda)$ and $(\tau, \xi) \stackrel{d}{=} \left(\frac{\eta}{\xi}, \xi\right)$, where $\eta \in Exp(1)$, and $\xi$ and $\eta$ are independent.
  \item[e)] For all $t > 0$, $P_{\xi}(x|\tau=t) = 0, \quad x \leq 0,$

  $P_{\xi}(x|\tau=t) =$
    $$ \frac{xt^2(\lambda + t)^3e^{-tx}(1-e^{-\lambda x}+\lambda x e^{-\lambda x})}{(\lambda + t)^3(1-e^{-at}-at \lambda e^{-at})+t^2(\lambda-t)(1-e^{-a(\lambda + t)})+ a (\lambda+t)t^3e^{-a(\lambda + t)}}, \,  x \in (0, a],$$
  $P_{\xi}(x|\tau=t) =$
   $$\frac{a\lambda xt^2(\lambda + t)^3 e^{-(\lambda+t)x}}{(\lambda + t)^3(1-e^{-at}-at \lambda e^{-at})+t^2(\lambda-t)(1-e^{-a(\lambda + t)})+ a (\lambda+t)t^3e^{-a(\lambda + t)}}, \,  x > a.$$
  \item[f)] The mean square regression function is $\mathbb{E}(\tau|\xi=x) = \frac{1}{x}$, $x > 0$.
   \item[g)] For $t > 0$, the mean square regression function is $\mathbb{E}(\xi|\tau=t) $
   $$= \frac{e^{-at}\left\{2e^{at}(1-6\lambda t^3)-(at+1)^2-1+e^{-a\lambda}t^3\left((a(\lambda+t)+1)^2+\lambda a^2(\lambda+t)-4a\lambda +1\right)\right\}}{t(\lambda + t)^3(1-e^{-at}-at \lambda e^{-at})+t^3(\lambda-t)(1-e^{-a(\lambda + t)})+ a (\lambda+t)t^4e^{-a(\lambda + t)}}.$$
\end{description}

{\bf Proof:} a) For $t > 0$,  by the integral form of the Total probability formula and (\ref{Def:1}) we obtain
\begin{eqnarray*}
  P_\tau(t) &=& \int_0^\infty P_\tau(t|\xi=x)P_\xi(x)dx = \int_0^a \frac{x}{a}e^{-xt}(1-e^{-\lambda x}+\lambda x e^{-\lambda x})dx + \int_a^\infty \lambda x e^{-(\lambda+t) x})dx\\
   &=& \frac{1}{at^2}\left(1-e^{-at}-at e^{-at}\right)+\frac{1}{a(\lambda + t)^2}\left(\frac{2\lambda}{\lambda+t}-1\right)(1-e^{-a(\lambda + t)}-a(\lambda + t)e^{-a(\lambda + t)})\\
   &+&\frac{\lambda}{(\lambda+t)^2}e^{-a(\lambda + t)}\\
   &=&  \frac{1}{at^2}\left(1-e^{-at}-at e^{-at}\right)+\frac{\lambda - t}{a(\lambda + t)^3}(1-e^{-a(\lambda + t)})+\frac{t}{(\lambda+t)^2}e^{-a(\lambda + t)}.
\end{eqnarray*}
Now, we compare it with (\ref{ExpMaxUExpCDF}) and complete the proof of a).

b) For $t > 0$, by the integral form of the Total probability formula  we obtain
\begin{eqnarray*}
  P_{\frac{\eta}{\xi}}(t) &=& \int_0^\infty P_{\frac{\eta}{\xi}}(t|\xi=x)P_\xi(x)dx = \int_0^\infty P_\frac{\eta}{x}(t)P_\xi(x)dx \\
  &=& \int_0^\infty xP_\eta(tx)P_\xi(x)dx =\int_0^a xe^{-xt}\frac{1}{a}(1-e^{-\lambda x}+ x \lambda e^{-\lambda x})dx  + \int_a^\infty x\lambda e^{-x(\lambda + t)}dx.
\end{eqnarray*}
The integrals are the same as in a), which means that for all $t > 0$, $P_\tau(t) = P_{\frac{\eta}{\xi}}(t)$. The rest follows by the uniqueness of the correspondence between p.d.f. and the probability law.

c) In order to obtain these moments we apply the Double expectation formula, and the formula for the moments of the exponential distribution.

d) follows by the formula $P_{\tau, \xi}(t,x) =  P_{\tau}(t|\xi = x)P_{\xi}(x)$, when we replace the p.d.f. of the exponential distribution, use its scaling property, and (\ref{MaxUExpDensity}).

e) can be proved by the Bayes' formula for the densities, (\ref{MaxUExpDensity}), (\ref{ExpMaxUExpDensity}),
and the p.d.f. of the Exponential distribution.

f) follows by the expectation of the Exponential distribution.

g) follows by the formula for the expectation, and e).
 \hfill $\Box$

 {\bf Definition 4.} \label{Def:4} We say that the rv. $(\tau_1, \tau_2, ..., \tau_k)$ has {\bf Multivatiate Exp-Max-U-Exp distribution of $II^{-nd}$ kind with parameters $a > 0$, and $\lambda > 0$}, if it has a joint p.d.f.

$P_{\tau_1, \tau_2, \ldots, \tau_k}(t_1, t_2, \ldots, t_k)$
$$ = \frac{\gamma(k+1,a(t_1 + t_2 + \ldots + t_k))}{a(t_1 + \ldots + t_k)^{k+1}} + \frac{\gamma(k+1, a(t_1 + \ldots + t_k + \lambda))}{a(t_1 + \ldots + t_k + \lambda)^{k+2}}(\lambda k - t_1 - \ldots - t_k)$$
$$ +\lambda k \frac{\Gamma(k, a(t_1 + \ldots + t_k + \lambda))}{(\lambda + t_1 + \ldots t_k)^{k+1}}, \quad t_1 > 0, t_2 > 0, \ldots, t_k > 0,$$
and $P_{\tau_1, \tau_2, \ldots, \tau_k}(t_1, t_2, \ldots, t_k) = 0$,  otherwise.

Briefly we will denote this in this way  $(\tau_1, \tau_2, \ldots, \tau_k) \in Exp-Max-U-Exp-II(a, \lambda)$.

{\bf Definition 5.} \label{Def:5} We say that the r.v. $T_n$ is {\bf{Erlang-Max-U-Exp distributed with parameters $n \in \mathbb{N}$, $a > 0$, and $\lambda > 0$}}, if it has a p.d.f.
$$P_{T_n}(t) = \frac{t^{n-1}}{a(n-1)!}\left(\frac{\gamma(n+1,at)}{t^{n+1}}+\frac{\gamma(n+1,a(\lambda+t))}{(\lambda+t)^{n+2}}(\lambda n - t)) + \lambda n a\frac{\Gamma(n,a(\lambda + t))}{(\lambda + t)^{n+1}}\right), $$
when $t > 0$, and $P_{T_n}(t) = 0$, otherwise. Briefly, we will denote this in this way  $T_n \in Erlang-Max-U-Exp(n; a, \lambda)$.

{\bf Theorem 3.} For $a > 0$, and $\lambda > 0$, if $\xi \in Max-U-Exp(a; \lambda)$ and for $x > 0$, $(\tau_1, \tau_2, ..., \tau_k|\xi=x)$ are independent identically $Exp(x)$ distributed r.vs., then,

\begin{description}
  \item[a)] $(\tau_1, \tau_2, \ldots, \tau_k) \in Exp-Max-U-Exp-II(a, \lambda)$.
  \item[b)] For $i = 1, 2, ..., k$, $\tau_i \in  Exp-Max-U-Exp-(a, \lambda)$.
  \item[c)] $(\tau_1, \tau_2, \ldots, \tau_k)  \stackrel{d}{=} \left(\frac{\eta_1}{\xi}, \frac{\eta_2}{\xi}, \ldots, \frac{\eta_k}{\xi}\right)$, where $\eta_1, \eta_2, \ldots, \eta_k$ are independent identically distributed (i.i.d.) $Exp(1)$, and independent on $\xi$.
 \item[d)]  $T_n := \tau_1 + \ldots + \tau_n \in Erlang-Max-U-Exp(n; a, \lambda)$. $T_n \stackrel{d}{=} \frac{\eta_1 + \eta_2 + \ldots + \eta_n}{\xi}$, where $\eta_1, \eta_2, \ldots, \eta_n$ are i.i.d. $Exp(1)$, and independent on $\xi$. $T_n \stackrel{d}{=} \frac{\theta_n}{\xi}$, where $\theta_n \in Gamma(n, 1)$ is independent on $\xi$.
   \item[e)]  For $p \in (0, 1)$,

 $\mathbb{E}(T_n^p)$
 $$  = \frac{\Gamma(p+n)}{(n-1)!\lambda^p}\left\{\frac{1}{a^p(1-p)}+\frac{\lambda^{p-1}}{a}\left((p+a\lambda)\Gamma(1-p,\lambda a)+(\lambda a)^{1-p}e^{-\lambda a} -p\Gamma(1-p)\right)\right\}.$$
  For $p \geq 1$, $\mathbb{E}(T_n^p) = \infty$.
  \end{description}

{\bf Proof:} a) For $t_1 > 0, t_2 > 0, \ldots, t_k > 0$ the integral form of the Total probability formula, and (\ref{MaxUExpDensity}) entail

$P_{\tau_1, \tau_2, \ldots, \tau_k}(t_1, t_2, \ldots, t_k)$
\begin{eqnarray*}
 &=& \int_0^\infty P_{\tau_1, \tau_2, \ldots, \tau_k}(t_1, t_2, \ldots, t_k|\xi = x)P_{\xi}(x)dx\\ &=& \frac{1}{a}\int_0^a x^k e^{-x(t_1 + t_2 + \ldots + t_k)} (1-e^{-\lambda x}+ x \lambda e^{-\lambda x}) dx +  \lambda \int_a^\infty x^k e^{-x(t_1 + t_2 + \ldots + t_k + \lambda)} dx\\
&=& \frac{\gamma(k+1,a(t_1 + t_2 + \ldots + t_k))}{a(t_1 + \ldots + t_k)^{k+1}} + \frac{\gamma(k+1, a(t_1 + \ldots + t_k + \lambda))}{a(t_1 + \ldots + t_k + \lambda)^{k+2}}(\lambda k - t_1 - \ldots - t_k)\\
& +&\lambda k \frac{\Gamma(k, a(t_1 + \ldots + t_k + \lambda))}{(\lambda + t_1 + \ldots t_k)^{k+1}}
\end{eqnarray*}
Otherwise  $P_{\tau_1, \tau_2, \ldots, \tau_k}(t_1, t_2, \ldots, t_k) = 0$. Now, we compare the last expression with Definition 4 and complete the proof of this point.

b) By condition $(\tau_1, \tau_2, ..., \tau_k|\xi=x)$ are i.i.d., therefore, for any fixed $i = 1, 2, ..., k$ we just can apply Theorem 1, a) and obtain immediately that $\tau_i \in  Exp-Max-U-Exp-(a, \lambda)$.

c) Consider $t_1 > 0, t_2 > 0, \ldots, t_k > 0$. analogously to the proof of a) we obtain the same expression as in a),

$P_{\frac{\eta_1}{\xi}, \frac{\eta_2}{\xi}, \ldots, \frac{\eta_k}{\xi}}(t_1, t_2, \ldots, t_k)$
\begin{eqnarray*}
 &=& \int_0^\infty P_{\frac{\eta_1}{\lambda}, \frac{\eta_2}{\lambda}, \ldots, \frac{\eta_k}{\lambda}}(t_1, t_2, \ldots, t_k|\xi = \lambda)P_{\xi}(\lambda)d\lambda \\
&=& \int_0^\infty P_{\eta_1, \eta_2, \ldots, \eta_k}(\lambda t_1, \lambda t_2, \ldots, \lambda t_k) \lambda^k P_{\xi}(\lambda)d\lambda \\
&=& \frac{1}{a}\int_0^a x^k e^{-x(t_1 + t_2 + \ldots + t_k)} (1-e^{-\lambda x}+ x \lambda e^{-\lambda x}) dx +  \lambda \int_a^\infty x^k e^{-x(t_1 + t_2 + \ldots + t_k + \lambda)} dx.
\end{eqnarray*}
Otherwise  $P_{\tau_1, \tau_2, \ldots, \tau_k}(t_1, t_2, \ldots, t_k) = 0$. The uniqueness of the correspondence between the p.d.f. and the probability distribution, together with Definition 4 complete the proof.

d) follows by the integral form of the Total probability formula, and the relation between the Erlang and Exponential distribution. The relation between Erlang, Exponential and Gamma distributions completes the proof.

e)  Consider $p \in (0, 1)$. By d) and the independence of $\eta_1, \eta_2, \ldots, \eta_k$ and $\xi$ we have

$\mathbb{E}(T_n^p)$
$$ = \mathbb{E}\left(\left(\frac{\eta_1 + \eta_2 + \ldots + \eta_n}{\xi}\right)^p\right) = \mathbb{E}((\eta_1 + \eta_2 + \ldots + \eta_n)^p) \mathbb{E}\left(\frac{1}{\xi^p}\right) = \frac{\Gamma(p+n)}{(n-1)!\lambda^p}\mathbb{E}\left(\frac{1}{\xi^p}\right),$$
where in the last equality we have used the well-known formula for the moments of $\eta_1 + \eta_2 + \ldots + \eta_n \in Gamma(n, \lambda)$.

Now, we use the definition for expectation together with (\ref{MaxUExpDensity}) and compute
$$\mathbb{E}\left(\frac{1}{\xi^p}\right) = \frac{1}{a^p(1-p)}+\frac{\lambda^{p-1}}{a}\left((p+a\lambda)\Gamma(1-p,\lambda a)+(\lambda a)^{1-p}e^{-\lambda a} -p\Gamma(1-p)\right),$$
which completes the proof.  \hfill $\Box$

\bigskip
\section{THE MIXED POISSON-MAX-U-EXP PROCESS}

{\bf Definition 6.} \label{Def:6} A r.v. $\theta$ has a {\bf{Mixed Poisson-Max-U-Exp distributed with parameters $a > 0$, and $\lambda > 0$}} if for $n = 0, 1, \ldots,$
\begin{equation}\label{MixedPoissonUEXP}
\mathbb{P}(\theta = n) = \frac{1}{n!}\left(\frac{\gamma(n+1,a)}{a} + \frac{\gamma(n+1,a(\lambda+1))}{a(\lambda+1)^{n+2}}(n\lambda-1)+ n\frac{\lambda\Gamma(n,a(\lambda+1))}{(\lambda+1)^{n+1}}\right).
\end{equation}
Briefly, $\theta \in MPMax-U-Exp(a, \lambda)$.

{\bf Definition 7.} \label{Def:7} Let $\mu(t): [0, \infty) \to [0, \infty)$ be a nonnegative, strictly increasing and continuous function, $\mu(0) = 0$, $\xi \in Max-U-Exp(a; \lambda)$ and $N_1$  be a Homogeneous Poisson process (HPP) with intensity $1$, independent on $\xi$. We call the random process
\begin{equation}\label{TheProcessN}
N := \{N(t), t\geq 0\}  = \{N_1(\xi \mu(t)), t \geq 0\}
\end{equation}
a {\bf{Mixed Poisson process with Max-U-Exp mixing variable}} or {\bf{MPMax-U-Exp process}}. Briefly $N \in MPMax-U-Exp(a, \lambda; \mu(t))$.

\medskip

{\bf Definition 8.} \label{Def:8} Let $n \in \mathbb{N}$. We say that a random vector $(N_1, N_2, \ldots, N_n)$ is {\bf Ordered Poisson-Max-U-Exp distributed with parameters $a > 0$, $\lambda > 0$, and $0 < \mu_1 < \mu_2 < ... < \mu_n$} if, for all integers $0 \leq k_1 \leq k_2 \leq \ldots \leq k_n$,

$\mathbb{P}(N_1 = k_1, N_2 = k_2, \ldots, N_n = k_n)$
$$ = \frac{\mu_1^{k_1}(\mu_2 - \mu_1)^{k_2 - k_1}\ldots(\mu_n - \mu_{n-1})^{k_n - k_{n - 1}}}{a k_1!(k_2 - k_1)!\ldots(k_n - k_{n-1})!}\left\{\frac{\gamma(k_n+1,a\mu_n)}{\mu_n^{k_n+1}}\right.$$
$$\left.+ \frac{\gamma(k_n+1,a(\lambda+\mu_n))}{(\lambda + \mu_n)^{k_n+2}}(\lambda k_n-\mu_n)+\lambda a k_n \frac{\Gamma(k_n,a(\lambda+\mu_n)) }{(\lambda+\mu_n)^{k_n+1}}\right\},$$
and $\mathbb{P}(N_1 = k_1, N_2 = k_2, \ldots, N_n = k_n) = 0$, otherwise. Briefly,
$$(N_1, N_2, \ldots, N_n) \in O_{PMUE}(a, \lambda; \mu_1, \mu_2, ..., \mu_n).$$

\medskip

{\bf Definition 9.} \label{Def:9} Let $n \in \mathbb{N}$. We say that a random vector $(N_1, N_2, \ldots, N_n)$ is {\bf Mixed Poisson-Max-U-Exp  distributed with parameters $a > 0$, $\lambda > 0$, and $0 < \mu_1 < \mu_2 < ... < \mu_n$} if, for all $m_1, m_2, \ldots, m_n \in \{0, 1, \ldots\}$,

$\mathbb{P}(N_1 = m_1, N_2 = m_2, \ldots, N_n = m_n) $
$$= \frac{\mu_1^{m_1}(\mu_2 - \mu_1)^{m_2}\ldots(\mu_n - \mu_{n-1})^{m_n}}{a m_1!m_2!\ldots m_n!}\left\{\frac{\gamma(m_1 + \ldots + m_n + 1,a\mu_n)}{\mu_n^{m_1 + \ldots + m_n +1}}\right.$$
$$+ \frac{\gamma(m_1 + \ldots + m_n+1,a(\lambda+\mu_n))}{(\lambda + \mu_n)^{m_1 + \ldots + m_n+2}}(\lambda (m_1 + \ldots + m_n)-\mu_n)$$
$$\left.+\lambda a (m_1 + \ldots + m_n) \frac{\Gamma(m_1 + \ldots + m_n,a(\lambda+\mu_n)) }{(\lambda+\mu_n)^{m_1 + \ldots + m_n+1}}\right\},$$
and $\mathbb{P}(N_1 = m_1, N_2 = m_2, \ldots, N_n = m_n) = 0$, otherwise. Briefly,
$$(N_1, N_2, \ldots, N_n) \in M_{PMUE}(a, \lambda; \mu_1, \mu_2, ..., \mu_n).$$

\medskip

The next statements are analogous to corresponding one in \cite{SMladen} and \cite{LMJ}.

{\bf Proposition 3.} If $(N_1, N_2, \ldots, N_n) \in O_{PMUE}(a, \lambda; \mu_1, \mu_2, ..., \mu_n)$, then
$$(N_1, N_2 - N_1, \ldots, N_n - N_{n-1}) \in M_{PMUE}(a, \lambda; \mu_1, \mu_2, ..., \mu_n).$$

{\bf Proposition 4.}  If $(N_1, N_2, \ldots, N_n) \in M_{PMUE}(a, \lambda; \mu_1, \mu_2, ..., \mu_n)$, then
$$(N_1, N_1 + N_2, \ldots, N_1 + N_2 + \ldots + N_n) \in O_{PMUE}(a, \lambda; \mu_1, \mu_2, ..., \mu_n).$$

In the next theorem we investigate the main properties of MPMax-U-Exp process $N$, defined in (\ref{TheProcessN}).

{\bf Theorem 4.} Let $a > 0$, $\lambda > 0$, and  $\mu(t): [0, \infty) \to [0, \infty)$ be a nonnegative, strictly increasing and continuous function, and $\{N(t), t \geq 0\} \in MPMax-U-Exp(a, \lambda; \mu(t))$.

\begin{description}
\item[a)] For all $t > 0$, $N(t) \in MPMax-U-Exp(a\mu(t), \frac{\lambda}{\mu(t)})$.
\item[b)] These processes are over-dispersed,
$$\mathbb{E}N(t) = \mu(t)\frac{a}{2} + \frac{\mu(t)}{a \lambda^2}(1 - e^{-\lambda a}),$$
$$\mathbb{D}N(t) = \mu(t)\frac{a}{2} + \frac{\mu(t)}{a \lambda^2}(1 - e^{-\lambda a}) + \mu^2(t)\left\{\frac{a^2}{12} - \frac{1}{\lambda^2}(1 + e^{-\lambda a}) \right.$$
$$\left. + \frac{4}{a \lambda^3}(1 - e^{-\lambda a})  + \frac{1}{a^2 \lambda^4}(1-e^{-\lambda a})^2 \right\}.$$
\item[c)] The probability generating function (p.g.f.) of the time intersections is
$$\mathbb{E}(z^{N(t)}) = \frac{1}{a\mu(t)(1-z)}(1 - e^{-\lambda a})$$
$$ - \frac{1}{a(\lambda + \mu(t)(1-z))}(1 - e^{-(\lambda + \mu(t)(1-z)) a} - \lambda a  e^{-(\lambda + \mu(t)(1-z)) a})$$
$$ + \frac{\lambda}{a(\lambda +\mu(t)(1-z))^2}(1 - e^{-(\lambda + \mu(t)(1-z)) a} - (\lambda +\mu(t)(1-z))  e^{-(\lambda + \mu(t)(1-z)) a}), \quad |z| < 1.$$
\item[d)] For $t > 0$, and $n = 0, 1, \ldots$, $P_{\xi}(x|N(t) = n) = 0$, when $x \leq 0$,

 $P_{\xi}(x|N(t) = n)$
 $$ = \frac{(\mu(t))^{n+1}x^ne^{-\mu(t)x}(1-e^{-\lambda x}+ x \lambda e^{-\lambda x})}{\frac{\gamma(n+1,a\mu(t))}{\mu(t)} + \frac{\mu^n(t)\gamma(n+1,a(\lambda+\mu(t)))}{(\lambda+\mu(t))^{n+2}}(n\lambda-\mu(t))+ na\frac{\lambda \mu^n(t)\Gamma(n,a(\lambda+\mu(t))}{(\lambda+\mu(t))^{n+1}}}, \,\,  x \in (0, a], $$
 $P_{\xi}(x|N(t) = n)$
 $$ = \frac{a (\mu(t))^{n+1}x^n\lambda e^{-x(\mu(t)-\lambda)}}{\frac{\gamma(n+1,a\mu(t))}{\mu(t)} + \frac{\mu^n(t)\gamma(n+1,a(\lambda+\mu(t)))}{(\lambda+\mu(t))^{n+2}}(n\lambda-\mu(t))+ na\frac{\lambda \mu^n(t)\Gamma(n,a(\lambda+\mu(t))}{(\lambda+\mu(t))^{n+1}}},\,\, x > a.$$
\item[e)] For $t > 0$, and $n = 0, 1, \ldots$, the mean square regression is

$\mathbb{E}(\xi|N(t) = n)$
$$ = \frac{\frac{\gamma(n+2,a\mu(t))}{a \mu^{n+2}(t)} + \frac{\gamma(n+2,a(\mu(t)+\lambda))}{a (\mu(t) + \lambda)^{n+3}}(\lambda(n+1)-\mu(t))+ \frac{\lambda (n+1)}{(\lambda+\mu(t))^{n+2}}\Gamma(n+1,a(\lambda+\mu(t)))}{\frac{\gamma(n+1,a\mu(t))}{a \mu^{n+1}(t) } + \frac{\gamma(n+1,a(\mu(t)+\lambda))}{a (\mu(t) + \lambda)^{n+2}}(\lambda n-\mu(t))+ \frac{\lambda n}{(\lambda+\mu(t))^{n+1}}\Gamma(n,a(\lambda+\mu(t)))}.$$
\item[f)] For all $k = 0, 1, \ldots$, $...$,
$$\mathbb{E}[N(t)(N(t)-1)(N(t) - k + 1)] = \frac{(a\mu(t))^k}{k+1} + \frac{k\mu^k(t)}{a\lambda^{k+1}}\gamma(k+1, a\lambda) + \frac{k\mu^k(t)}{\lambda^k}\Gamma(k, \lambda a).$$
\item[g)] For all $n \in \mathbb{N}$, and $0 \leq t_1 \leq t_2 \leq \ldots \leq t_n,$
$$(N(t_1), N(t_2), \ldots, N(t_n)) \in O_{MPUE}(a, \lambda; \mu(t_1), \mu(t_2), ..., \mu(t_n)).$$
\item[h)] For all $n \in \mathbb{N}$, and $0 \leq t_1 \leq t_2 \leq \ldots \leq t_n$,
$$(N(t_1), N(t_2) - N(t_1), \ldots, N(t_n) - N(t_{n-1})) \in M_{MPUE}(a, \lambda; \mu(t_1), \mu(t_2), ..., \mu(t_n)).$$
\item[i)] Denote by $\tau_1, \tau_2, \ldots$ the inter-occurrence times of the counting process $N$. Then, $\tau_1, \tau_2, \ldots$ are dependent and $Exp-Max-U-Exp(a; \lambda)$ distributed.
\item[j)] For $n \in \mathbb{N}$, if $T_n$ is the moment of occurrence of the $n$-th event of the counting process $N$, then $T_n \in Erlang-Max-U-Exp(n; a, \lambda)$.
\end{description}

{\bf Proof:} a) Consider $t > 0$ and $n \in \mathbb{N}\cup \{0\}$. By Definition \ref{Def:7} we have that $\mathbb{P}(N(t) = n) = \mathbb{P}(N_1(\xi \lambda(t))= n)$.

The integral form of the Total probability formula, and the independence between the random process $N_1$, and the r.v. $\xi$ entail,
$$\mathbb{P}(N(t) = n) = \int_0^\infty \mathbb{P}(N_1(\xi \mu(t))= n|\xi=x)P_\xi(x) dx = \int_0^\infty \mathbb{P}(N_1(x \mu(t))= n)P_\xi(x)dx.$$

Now, by using the definition for Poisson distribution and the Definition \ref{MaxUExpDensity}, for $n = 0, 1, \ldots$ we obtain,

$\mathbb{P}(N(t) = n)$
\begin{eqnarray*}
 &=& \int_0^\infty \frac{(x\mu(t))^n}{n!}e^{-x\mu(t)}P_\xi(x) dx   \\
   &=& \int_0^a\frac{(x\mu(t))^n}{n!}e^{-x\mu(t)}\frac{1}{a}(1-e^{-\lambda x}+ x \lambda e^{-\lambda x}) dx + \int_a^\infty \frac{(x\mu(t))^n}{n!}e^{-x\mu(t)}\lambda e^{-\lambda x} dx \\
   &=& \frac{\gamma(n+1,a\mu(t))}{a \mu(t)n!} + \frac{\mu^n(t)\gamma(n+1,a(\lambda+\mu(t)))}{a(\lambda+\mu(t))^{n+2}n!}(n\lambda-\mu(t))+ n\frac{\lambda \mu^n(t)\Gamma(n,a(\lambda+\mu(t))}{(\lambda+\mu(t))^{n+1}n!}.
\end{eqnarray*}
By definition \ref{Def:6}, the last expression is exactly the p.m.f. of $MPMax-U-Exp(a\mu(t), \frac{\lambda}{\mu(t)})$ distributed r.v. The rest follows by the uniqueness of the correspondence between the p.m.f. and the probability law of the r.v.

b) follows by the general formulae for the mean and the variance of Mixed Poisson distribution which could be seen for example in Proposition 2.1.i) and ii) in \cite{GrandelMixed} and Theorem 1 b) and c).

c) Let $z < 1$ and $t \geq 0$. By the Double expectation formula we have the general formula for the p.g.f. of a Mixed poisson process, which could be seen for example in  \cite{GrandelMixed} for the case when $\mu(t) \equiv t$, $t > 0$. It is
$\mathbb{E}(z^{N(t)}) = \mathbb{E}(e^{-\mu(t)\xi(1-z)}).$
Now, Theorem 1, d), completes the proof of this statement.

d) The Bayes rule, a), Definition \ref{Def:6}, (\ref{MaxUExpDensity}) and the definition for Poisson distribution entail the desired result.

e) In order to prove this statement we use the definition for the expectation and d).

f)  Remark 2.1, p. 15 in \cite{GrandelMixed} expresses the relation between these factorial moments and the moments of the mixing variable. Now, we use Theorem 1, b) and complete the proof of this point.

g) and h) are analogous to the proves of the analogous results in \cite{SMladen} and \cite{LMJ}.

i) follows by Definition 4, the properties of the HPP, and Theorem 2, a).

j) follows by Definition 5, the properties of the HPP, and Theorem 2, d).  \hfill $\Box$

{\it{Notes:}} 1. As it is noticed in Grandel \cite{GrandelMixed}, in the case $\mu(t) = t$, $t > 0$, any Mixed Poisson process is a birth process with transition intensities given by $\mathbb{E}(\xi|N(t) = n)$. The converse is not true.

 In the general case for $\mu = \mu(t)$, $t > 0$, the transition intensities are analogous, however, first we need to apply the non-random time-change $\mu$, to the initial birth process.

2.  As far as any Mixed Poisson process with $\mu(t) = t$, $t > 0$, for any $0 < s < t$ has Binomial conditional distributions $(N(s)|N(t) = n) \in Bi\left(n, \frac{s}{t}\right)$, (see Grandel \cite{GrandelMixed}), p. 98, in the general case for $\mu$ we have
$$(N(s)|N(t) = n) \in Bi\left(n, \frac{\mu(s)}{\mu(t)}\right).$$

\bigskip
\section{CONCLUSIONS}
Mixed Poisson processes represent a generalization of homogeneous ones, allowing the rate $\lambda$ to be a random variable. The distribution $\lambda$  is called the structure distribution and may be regarded as a prior distribution. This work considers a new structure distribution, called Max-U-Exp distribution. It is the distribution of the maximum of two random variables - Uniform and Exponential. Properties of this distribution are considered and an algorithm for estimating its parameters is developed. It is based on a combination of the method of moments and least square method. Exp-Max-U-Exp distribution is defined and thoroughly investigated. It arises in a natural way as a distribution of the inter-arrival times in the Mixed Poisson process with Max-U-Exp mixing variable. The distribution of the renewal moments (arrival times) is called Erlang-Max-U-Exp and is defined via its probability density function.
By using mainly the previous general theory about Mixed Poisson processes, developed  by Grandel \cite{GrandelMixed}, and Karlis and Xekalaki \cite{KarlisandXekalaki}, and analogously to Jordanova et al. \cite{SMladen}, and Jordanova and Stehlik \cite{LMJ} we define and investigate the properties of the new random vectors and random variables, which are related with this particular case of a Mixed Poisson process. The paper shows new explicit relations between the considered random elements. In an analogous way, many different univariate and multivariate distributions could be defined, and different relations between the new classes of probability laws could be explained.

\bigskip
\section{ACKNOWLEDGMENTS} The work was supported by the Scientific Research Fund in Konstantin Preslavsky University of Shumen, Bulgaria under Grant Number RD-08-35/18.01.2023 and project Number 2023 - FNSE – 04, financed by Scientific Research Fund of Ruse University.

An improved version of this article has been accepted by AIP Conference Proceedings, 49th International Conference Applications of Mathematics in Engineering and Economics, 10 - 16 June 2023, Sozopol, Bulgaria.

\bibliographystyle{plain}

\end{document}